\makeatletter
\@namedef{subjclassname@2020}{\textup{2020} Mathematics Subject Classification}
\makeatother

\def\bbuildrel#1_#2^#3{\mathrel{\mathop{\kern 0pt#1}\limits_{#2}^{#3}}}

\def\NN{\mathbb N}

\def\RR{\mathbb R}

\def\bs{\bigskip}
\def\ms{\medskip}
\def\ss{\smallskip}

\def\w{\thinspace\hbox{\hsize 14pt \rightarrowfill}\thinspace}

\def\0{\hbox{$\emptyset$}}
\def\A{\hbox{$\mathcal A$}}

\def\s{\hbox{$\sigma$}}

\def\sub{\subseteq}
\def\G{\hbox{$\mathcal G$}}

\def\N{\hbox{$\mathcal N$}}
\def\F{\hbox{$\mathcal F$}}

\def\B{\mathscr{B}}

\def\D{\mathscr{D}}

\def\E{\mathscr{E}}
\def\H{\mathscr{H}}
\def\G{\mathscr{G}}

\def\U{\mathscr{U}}
\def\V{\mathscr{V}}
\def\W{\mathscr{W}}

\documentclass[a4paper,12pt]{amsart}

\usepackage{amssymb,enumerate,mathrsfs}
\usepackage[utf8]{inputenc} 

\usepackage{enumerate}

\theoremstyle{plain}

\newcommand{\co}{\mathfrak c}
\newcommand{\be}{\mathfrak b}
\newcommand{\de}{\mathfrak d}

\newcommand{\cantor}{2^{\NN}}
\newcommand{\ctc}{2^{\NN}\times 2^{\NN}}

\newcommand{\baire}{\NN^{\NN}}
\newtheorem{theorem}{Theorem}[section]

\newtheorem{lemma}[theorem]{Lemma}

\newtheorem{proposition}[theorem]{Proposition}

\numberwithin{equation}{section}

\begin{document}

	\title{
    On SierpiŃski sets, Hurewicz spaces and Hilgers functions}

	\author{Witold Marciszewski, Roman Pol and Piotr Zakrzewski}
	
	
	\address{Institute of Mathematics, University of Warsaw,  Banacha 2, 02-097 Warsaw, Poland}
	\email{wmarcisz@mimuw.edu.pl}
	\email{pol@mimuw.edu.pl}
	\email{piotrzak@mimuw.edu.pl}

	
	
	\date{\today}

	\begin{abstract} 
		The Hurewicz property is a classical generalization of $\s$-compactness and Sierpiński sets (whose existence follows from CH) are standard examples of non-\s-compact Hurewicz spaces. We show, solving a problem stated by Szewczak and Tsaban in \cite{sz-ts-1}, that for each Sierpiński set $S$ of cardinality at least $\be$ there is a Hurewicz space $H$ with $S\times H$ not Hurewicz.

        Some other questions in the literature concerning this topic are also answered.
	\end{abstract}

	\maketitle
	
	\section{Introduction}\label{intro}

    The Hurewicz property, introduced in \cite{hur-2}, is a classical generalization of \s-compactness which attracted much attention in point-set topology and set theory (see \cite{ts} for an excellent self-contained introduction to this topic).

    Sierpiński sets (whose existence follows from CH) provided first examples of non-\s-compact sests with the Hurewicz property.

    In this article we present solutions to some open problems from the literature of the subject. In particular:
    \begin{itemize}
\item Solving Problem 7.5 formulated by Szewczak and Tsaban in \cite{sz-ts-1} we prove  (cf. Theorem \ref{IMG-2, 27-12-2024}) that for every Sierpiński set $S\sub \cantor$  of cardinality at least $\be$ there is a Hurewicz space $H\sub\cantor$ such that $S\times H$ is not a Hurewicz space.

\item  We provide a solution of the second part of Problem 3 stated by Banakh and Zdomskyy in \cite{Ba-Zd} by showing (cf. Theorem \ref{IMG, 9-12-2024}) that, assuming $V=L$, there is a set in $\cantor$ with the Borel-separation property but without the Analytic-separation property. 

\item  We give  a negative answer (cf. Theorem  \ref{RP, 23-11-2024}) to Problem 6.13 formulated by Sakai and Scheepers in \cite{sak-sch} by proving that, assuming CH, for each Sierpiński set $S\sub \cantor$ there is a Hurewicz $C$-space $H\sub I^\NN$ such that $S\times H$ is not a $C$-space.

\item Solving Problem 6.6 formulated by Sakai and Scheepers in \cite{sak-sch} we confirm their conjecture by showing (cf. Theorem  \ref{IMG, 27-11-2024}) that,   assuming $\be=\co$, there is a $C$-space $E\sub I^\NN$ such that $E^n$ is a Hurewicz space for $n=1,2,\ldots$ but $E\times \NN^\NN$ is not a $C$-space.
	
\end{itemize}

Some key ideas of our approach go back to the seminal paper by E. Michael \cite{mich}. These ideas were used in the literature in connection with the Baire category, involving Lusin sets and Menger spaces, cf. the references in \cite{p-p}.

The Hilgers functions, vital in some of our reasonings, are recalled in subsection \ref{Hilgers functions}.
    
	\section{Terminology and some auxiliary results}\label{notation}
	
	\subsection{Notation}

We identify the Cantor set with the product $\cantor$. The constant zero sequence in $\cantor$ is denoted by $\bar{0}$   and $\bar{Q}$ is the subspace of $\cantor$ consisting of all eventually zero sequences in $\cantor$ - a homeomorphic copy of the space of rationals $\mathbb{Q}$. Then $\bar{P}$ denotes $\cantor \setminus \bar{Q}$ - a homeomorphic copy of the Baire space $\baire$. 

For a subset $A$ of a product $X \times Y$, and $x\in X, y\in Y$, $A_x = \{y: (x,y)\in A\}$ is the vertical section of $A$ at $x$, and the horizontal section $\{x: (x,y)\in A\}$ of $A$ at $y$ is denoted by $A^y$.   

The smallest cardinality of a subset of $\cantor$ which is  nonmeasurable with respect to the standard probability product measure $\lambda$ on $\cantor$ is denoted by non($\N$) and {\rm cof}($\N)$ denotes the smallest cardinality of a base of the \s-ideal of all $\lambda$-null sets (shortly: null sets). The smallest cardinality of a covering on $\cantor$ by null sets is denoted by ${\rm cov}(\N)$.

The smallest cardinality of a  subset of $\baire$ which is unbounded (dominating, respectively) in the ordering  $\leq^*$  of eventual domination is denoted by  $\mathfrak{b}$ ($\mathfrak{d}$, respectively), cf.  \cite{Bla}.
    
	\subsection{Sierpiński sets
    }\label{subs_Sierp}

Let us recall that $S\sub \cantor$ is a {\sl Sierpiński set}, if it is uncountable and has countable intersection with every null set in $\cantor$, cf. \cite{mi-1}. More generally, for an uncountable cardinal $\kappa$ we call $S\sub \cantor$ a {\sl $\kappa$-Sierpiński set}, if $|S|\geq\kappa$  and $|S\cap N|<\kappa$ for  every null set $N$ in $\cantor$, cf. \cite{Mi-Ts-Zd}. In particular, a $\kappa$-Sierpiński set exists  under the assumption that ${\rm cov}(\N) ={\rm cof}(\N)=\kappa$.

\ss

We shall use the following observations which presumably belong to a folklore.

	\begin{lemma}\label{Sierp in B}
		Let $S$ be a Sierpiński set in $\cantor$. For any Borel set $B\subseteq \cantor$ of positive measure there exists a Sierpiński set $E\subseteq B\setminus S$.
	\end{lemma}
	
	\begin{proof} 
  
Let us recall  that there is a Borel isomorphism $\varphi$ of $\cantor$ onto $B$, preserving $\lambda$-null sets, cf. \cite[Theorem 17.41]{ke}. It follows that if  $B\cap S$ is  countable, then $E=\varphi(S)\setminus S$ is a required Sierpiński set. If, on the other hand, $B\cap S$ is uncountable, then $S_1=\varphi^{-1}(B\cap S)$ is a Sierpiński set and if $T$ is any Sierpiński set in $\cantor$ disjoint from $S_1$, then $E=\varphi(S_1)$ satisfies our requirements. 

 \ss 
    
	To check that, indeed, there is a Sierpiński set $T$ in $\cantor$ disjoint from $S$,  we shall split the argument into two cases according to the 
    validity or non-validity of CH.
	
	\ss

	{\bf Case (A):}  $2^{\aleph_0}=\aleph_1$.

    \ss 

        Let $N_\alpha$, $\alpha<\omega_1$ be  $\lambda$-null $G_\delta$-sets in $\cantor$ such that every $\lambda$-null set in $\cantor$ is a subset of some $N_{\alpha}$. We inductively choose 
        $$
        e_{\alpha}\in \cantor\setminus  (S\cup \bigcup_{\beta<\alpha} N_{\beta}\cup \{e_\beta:\beta<\alpha\}),
        $$
        and then $T=\{e_{\alpha}:\alpha<\omega_1\}$ is a required Sierpiński set.
        
\ss

	{\bf Case (B):}  $2^{\aleph_0}>\aleph_1$.

    \ss 

    Let us take $M\sub S$ with $|M|=\aleph_1$. We claim that there is $x\in\cantor$ such that if we let $T=x+M$, then $T\cap S=\emptyset$ and  $T$ is a required Sierpiński set. Indeed, otherwise $S\cap (x+M)\neq \emptyset$ for each $x\in\cantor$ or, equivalently, $M-S=\bigcup_{x\in M}(x-S)=\cantor$. But this is impossible, since the intersection of $M-S$ with any null set in $\cantor$ has cardinality at most $\aleph_1$.
        \end{proof}
	
	\begin{lemma}\label{Sierp full measure}
		Let $S$ be a Sierpiński set in $\cantor$. Then there exists a Sierpiński set $T\subseteq \cantor\setminus S$ of full outer measure $\lambda^*(T) = 1$.
	\end{lemma}
	
	\begin{proof}
		Let $\E$ be a maximal collection of pairs $(E,E^*)$, where $E\sub \cantor\setminus S$ is a Sierpiński set (cf. Lemma \ref{Sierp in B}), $E^*$ is a Borel set with $E\sub E^*$ and $\lambda(E^*)=\lambda^{*}(E)$, and for any distinct  $(E_1,E_1^*)$, $(E_2,E_2^*)$ in $\E$,  $E_1^*\cap E_1^*=\emptyset$. The family $\E$ is countable, and by Lemma \ref{Sierp in B}, if we let $T=\bigcup\{E: (E,E^*)\in \E \}$, then $T$ is a required Sierpiński set.
	\end{proof}

	\subsection{Hurewicz and Menger spaces}\label{Hur_and_Meng} 
    Let us recall that  a separable metrizable space $X$  is {\it a  Hurewicz space}, if for each sequence $\U_1, \U_2,\ldots$ of open  covers of $X$,
there are finite subfamilies $\F_n\sub \U_n$ such that $X = \bigcup_n \bigcap_{m\geq n}(\bigcup \F_m)$. 
By a theorem of Hurewicz (cf. \cite{hur-2}) this is equivalent  
 to the statement that for every continuous function $f:X\w \RR^\NN$ the image of $X$ is bounded (in the sense that the subset of $\baire$ consisting of sequences of the form $(\lceil |f(x)(n)|\rceil)_{n\in\NN}$
 is bounded in the ordering  $\leq^*$  of eventual domination). Any \s-compact space  is a Hurewicz space and there exist (in ZFC)  Hurewicz, non-\s-compact spaces (cf. subsection \ref{bscale}). 
 A useful application of the above characterization is the fact that a separable metrizable space $X$ which is a union of less than $\be$ Hurewicz subspaces $X_\alpha$, is Hurewicz, since its image under any continuous function $f:X\to \RR^\NN$ is bounded as the union of less than $\be$ bounded sets  of the form $f(X_\alpha)$; in particular, if $|X|<\be$, then $X$ is a Hurewicz space.

We shall often apply the following characterization of Hurewicz spaces  (cf. \cite[Theorem 5.7]{j-m-s-s}): {\sl a subspace $X$ of a compact, metrizable space $K$ is a Hurewicz space if and only if for any $G_\delta$-set $G$ in $K$ containing $X$, there is an $\sigma$-compact set $F$ in $K$ such that $X\sub F\sub G$}.

The characterization yields readily that any
$\be$-Sierpiński set $S$ in $\cantor$ is  a non-\s-compact Hurewicz space (note, however, that the existence of $\be$-Sierpiński sets cannot be proved in ZFC). Indeed, if  $G$ is a $G_\delta$-set  in $\cantor$ containing $S$ and $F_1\sub G$ is a \s-compact set of measure $\lambda(G)$, then since $S$ is a $\be$-Sierpiński set, we have that $|S\setminus F_1|<\be$. Consequently, $S\setminus F_1$ is a Hurewicz space so it can be covered by a \s-compact set $F_2\sub G$ and if we let $F=F_1\cup F_2$, then $F$ is a \s-compact set with $S\sub F\sub G$ which witnesses that $S$ is a Hurewicz space. 
\ss

 Let us also recall that a separable metrizable space $X$ is {\it a Menger space}, if for each sequence $\U_1, \U_2,\ldots$ of open  covers of $X$,
there are finite subfamilies $\F_n\sub \U_n$ such that $X = \bigcup_n (\bigcup \F_n)$ (cf. \cite{meng}, \cite{hur-1}).  Clearly, if $X$ is \s-compact, then it is a Menger space and every Hurewicz space is also a Menger space
but there are (in ZFC) Menger spaces which are not Hurewicz, cf \cite{ts}.
\ss 

\subsection{$\be$-scale set $\mathfrak{B}$ of Bartoszyński and Shelah} \label{bscale} 
In \cite{b-s} Bartoszyński and Shelah gave a ZFC example of a Hurewicz, non-\s-compact subspace $\mathfrak{B}$ of $\cantor$ with the following properties:
\begin{enumerate}[{(BS1)}]\label{bscalep} 
\item\label{bscale1} $\mathfrak{B}$ contains $\bar{Q}$ and has cardinality $|\mathfrak{B}| = \be$;
\item\label{bscale2} for each $\sigma$-compact subset $K$ of $\bar{P}$, we have $|\mathfrak{B}\cap K| < \be$ and $\mathfrak{B} \setminus K$ is a Hurewicz space.
\end{enumerate}

More precisely, $\mathfrak{B}$  is the union of $\bar{Q}$ with a copy $A$ (under a homeomorphism  between $\NN^{\NN}$  and  $\bar{P}$) of a {\sl $\be$-scale}, i.e.,  any well-ordered by eventual domination and unbounded subset of $\baire$. 

Following \cite{Mi-Ts-Zd} we shall call $\mathfrak{B}$ {\sl a $\be$-scale set} (clearly, 
$\mathfrak{B} \setminus K$ is also a $\be$-scale set, for any \s-compact $K\sub \bar{P}$).

A simpler proof that every $\be$-scale set is a Hurewicz space can be found in \cite[Remark 4.2]{p-z}, where the argument is based on the following 
observation concerning the set $\mathfrak{B}$ ((BS3) below implies (BS2), cf. \cite[Remark 4.2]{p-z}, and the reverse implication follows readily from the characterization in subsection \ref{Hur_and_Meng}):
\begin{enumerate}[{(BS3)}]
\item\label{bscale3} for each $\sigma$-compact subset $K$ of $\bar{P}$, there exists a $\sigma$-compact subset $L$ of $\cantor$ such that $\bar{Q}\subseteq L,\ L\cap K = \emptyset$, and $|\mathfrak{B} \setminus L| < \be$. 
\end{enumerate}	

Note that condition (BS2) implies that the subspace $A = \mathfrak{B}\setminus \bar{Q}$ is not Hurewicz (cf. the characterization of Hurewicz spaces from subsection \ref{Hur_and_Meng}).
	
If $\be=\de$, then we can require that the set $A$ used in the construction of space $\mathfrak{B}$ is a copy of a scale, i.e.,  a well-ordered by eventual domination, cofinal and unbounded subset of $\baire$. In this case, the subspace $A$ of $\mathfrak{B}$ is not a Menger space, cf. \cite[Lemma 1.4]{ts}.

	\subsection{$\lambda$-spaces and $\lambda'$-sets}
    
    Let us recall that a subspace $X$ of a compact metrizable space $K$ is {\sl a $\lambda$-space} if  every countable set $N\sub X$ is relatively $G_{\delta}$ in $X$ and is {\sl  a  $\lambda'$-set in $K$} if every countable set $N\sub K$ is relatively $G_{\delta}$ in $X \cup N$.

\ss

Clearly, if $X$ is a  $\lambda'$-set in $K$, then it is a $\lambda$-space but the opposite implication is false (cf.\ \cite[Theorem 5.6]{mi-1}). 
Let us recall the following two observations.

\begin{lemma}\label{lambda'}
If $X$ is a Hurewicz $\lambda$-space contained in a compact metrizable space $K$, then $X$ is a  $\lambda'$-set in $K$.
\end{lemma}

    \begin{proof}
Let as assume that $X$ is a Hurewicz $\lambda$-space contained in $K$ and let $N$ be an arbitrary countable subset of $K$. Let $N_1=N\cap X$ and $N_2=N\setminus X$. 

Since $X$ is a $\lambda$-space, there is a $G_\delta$-set  $G_1$ in $K$ such that $N_1=G_1\cap X$.

Since $X$ is a Hurewicz space and $K\setminus N_2$ is a $G_\delta$-set in $K$ containing $X$, there is a $G_\delta$ set $G_2$  in $K$ containing $N_2$ and disjoint from $X$.

Let $G=G_1\cup G_2$. Then $G$ is a $G_\delta$-set in $K$ and $N=G\cap (X\cup N)$.
 \end{proof}

\begin{lemma}\label{Sierp_lambda'}
If $S$ is a Sierpiński set in $\cantor$, then $S$ is a $\lambda$-space and hence a $\lambda'$-set in $\cantor$.
\end{lemma}

   \begin{proof}
Let $N$ be a countable subset of $S$, and let $G$ be a $G_\delta$-set in $\cantor$ covering $N$, with  $\lambda(G) =0$. Then $G\cap S$ is a countable $G_\delta$-set in $S$ containing $N$, which readily implies that $N$ is a $G_\delta$-set in $S$.
 \end{proof}
	
	\subsection{Hilgers functions}\label{Hilgers functions} In this note we shall frequently use the following diagonal construction, going back to Hilgers \cite{hil}.

    Let $S$ be a subset of a set $X$, $A\sub X\times Y$ be a subset of the product $X\times Y$ of sets $X,\ Y$ such that the projection $\pi_X(A)$ contains $S$, and let $\mathcal{F} = \{F_\alpha: \alpha < \co\}$ be a family of subsets of  $X\times Y$ with $S\sub\pi_X(F_\alpha)$ for $\alpha<\co$.
    Given a partition $\mathcal{P} = \{S_\alpha: \alpha < \co\}$ of the set $S$ into non-empty sets such that $S_\alpha \cap S_\beta = \emptyset$ for $\alpha\neq \beta$, we define a function $f: S\to Y$ in the following way: given  $x\in S_\alpha$, where $\alpha < \co$, we pick $f(x)\in A_x\setminus (F_\alpha)_x$, whenever such choice is possible, and we choose  $f(x)\in A_x$ arbitrarily, if $A_x\sub (F_\alpha)_x$.  
	
    
    We shall say, cf. \cite{p-p}, that $f$ is a  
	{\sl Hilgers function} associated with the set $A$, the family $\mathcal{F}$ and the partition  $\mathcal{P}$. If the partition  $\mathcal{P}$ consists of singletons, we simply say that $f$ is associated with  $A$ and $\mathcal{F}$. 
    
    One can easily verify that  $Gr(f)$, the graph of $f$, has the following property
	\begin{equation*}\label{Hilg_funct}\tag{H1}
		\mbox{for any } \alpha < \co, \mbox{ if } Gr(f)\subseteq F_\alpha, \mbox{ then }  A_x\subseteq (F_\alpha)_x   \mbox{ for all } x\in S_\alpha.
	\end{equation*}

	\section{No Sierpiński set of cardinality at least $\be$ is productively Hurewicz}
	\bs  
	
	The following result provides a positive answer (in a strong form) to Problem 7.5 in \cite{sz-ts-1} (repeated as Problem 5.5 in \cite{sz-ts-2}). 
	
	\begin{theorem}\label{IMG-2, 27-12-2024}
		For every Sierpiński set $S\sub \cantor$ of cardinality at least $\be$ there is a Hurewicz $\lambda'$-set $H$ in $\cantor$ such that $S\times H$ is not Hurewicz; if $\be=\de$, one can have $S\times H$ not Menger.
	\end{theorem}
	
	\begin{proof}
Let $T$ be a Sierpiński set in $\cantor$ disjoint from $S$, with $\lambda^*(T) = 1$, given by Lemma \ref{Sierp full measure} (recall that $\lambda^*$ is the outer measure on $\cantor$).

We shall follow closely the reasoning from 	\cite[Example 4.1 and Remark 4.2]{p-z} concerning the $\be$-scale set $\mathfrak{B}$ described in 
Subsection \ref{bscale}.

We put $A = \mathfrak{B}\setminus \bar{Q}$. Let $S'$ be a proper subset of $S$ of cardinality $\be$, and let $g: S\to A$ 
map $S\setminus S'$ to a point $a\in A$, and $S'$ bijectively onto $A\setminus \{a\}$.

We will prove that the set
\begin{equation}\label{3.1.1}
H = Gr(g)\cup (T\times \bar{Q}) \subseteq (S\cup T)\times \mathfrak{B}
\end{equation}
has the required properties (identifying $\cantor \times \cantor$ with $\cantor$, we can treat $H$ as a subset of $\cantor$).

Let $\pi_i: \cantor \times \cantor \to \cantor,\ i=1,2$, denote the projections onto the first and second axis, respectively.

First, we will verify that 
\begin{equation}\label{3.1.2}
H \mbox{ is a Hurewicz space.}
\end{equation}
To that end, fix a $G_\delta$-set $G$ in $\cantor \times \cantor$ containing $H$. We have to find a $\sigma$-compact set lying between $H$ and $G$ (cf. subsection \ref{Hur_and_Meng}).

Since $\bigcap\{G^q: q\in \bar{Q}\}$ is a $G_\delta$-set containing $T$, and 
$T$ being a Sierpiński set, is a Hurewicz space
(cf.\ subsection \ref{Hur_and_Meng}), there is a $\sigma$-compact set $F$ with $T \subseteq F \subseteq \bigcap\{G^q: q\in \bar{Q}\}$. Therefore, for a $\sigma$-compact set $F\times \bar{Q}$ we have
\begin{equation}\label{3.1.3}
T\times \bar{Q} \subseteq F\times \bar{Q} \subseteq G.
\end{equation}
Let 
\begin{equation}\label{3.1.4}
K = \pi_2((F\times \cantor)\setminus G).
\end{equation}
Then $K$ is a $\sigma$-compact set in $\bar{P} = \cantor \setminus \bar{Q}$, cf. (\ref{3.1.3}). By property (BS3) of $\mathfrak{B}$ (cf.\ subsection \ref{bscale}) we can find a $\sigma$-compact subset $L$ of $\cantor$ such that
\begin{equation}\label{3.1.5}
\bar{Q}\subseteq L \subset (\cantor\setminus K)\mbox{ and } |\mathfrak{B} \setminus L| < \be.
\end{equation}
 By (\ref{3.1.3}) and (\ref{3.1.4}) we have
\begin{equation}\label{3.1.6}
T\times \bar{Q} \subseteq F\times L \subseteq G.
\end{equation}
Since $\lambda^*(T) = 1$ and $\cantor\setminus F$ is disjoint from $T$, we have $\lambda(\cantor\setminus F) = 0$, 
and $S$ being a Sierpiński set, $|(\cantor\setminus F)\cap S| \le \aleph_0$. By the definition of $H$ (cf.\ (\ref{3.1.1})), 
 \mbox{the set } $N=[(\cantor\setminus F)\times \cantor]\cap H$ \mbox{ is countable.}
Consequently, by (\ref{3.1.6}), $F_1=(F\times L)\cup N$ is a $\sigma$-compact subset of $G$ containing $H\cap (\cantor\times L)$.


Finally, letting $X=H\setminus (\cantor\times L)$, by  (\ref{3.1.1}) and (\ref{3.1.5}), we have $X \subseteq S\times(\mathfrak{B} \setminus L)$. Since $S$ is a Sierpiński set and $|\mathfrak{B} \setminus L| < \be$ (cf.\ (\ref{3.1.5})), $X$ is the union of less than $\be$ Hurewicz sets of the form $X^y\times\{y\}$, where $y\in \mathfrak{B}\setminus L$ and $X^y\sub S$ (cf. subsection \ref{Hur_and_Meng}). It follows that $X$, being a Hurewicz space,  is contained in a $\sigma$-compact subset $F_2$ of $G$.
Therefore,  $F_1 \cup F_2$  is a $\sigma$-compact subset of $G$ containing $H$, witnessing (\ref{3.1.2}).

\ss 
 
To prove that $H$ is a $\lambda$-space, it is enough to observe that $H$ is a set with countable vertical sections $H_x$, $x\in S\cup T$, over a Sierpiński set. Indeed,  if $N$ be a countable subset of $H$, then $\pi_1(N)$ is a countable subset of $S\cup T$, so because $S\cup T$ is a $\lambda$-space (cf.\ Lemma \ref{Sierp_lambda'}), there is a $G_\delta$-set $L$ in $\cantor$ with $L\cap (S\cup T)= \pi_1(N)$. Then $L \times \cantor$ is a $G_\delta$-set in $\cantor \times \cantor$ and $(L \times \cantor) \cap H$ is a countable $G_\delta$-set in $H$ containing $N$. 

Since we have already established that $H$ is a Hurewicz space, it is a $\lambda'$-set in $\cantor \times \cantor$, by Lemma \ref{lambda'}.

\ss 

It remains to check that
$S\times H$  is not a Hurewicz space.
But this follows readily from the fact that the product $S\times H$ contains a closed copy $\{(s,(x,y))\in S\times H: s= x\}$ of $Gr(g)$ and $\pi_2(Gr(g)) = g(S) = A$ is not Hurewicz (cf. the end of subsection \ref{bscale}). If $\be=\de$, then the set $A$ is not Menger (cf. \cite[Lemma 1.4]{ts}), implying that the product $S\times H$ is also not Menger.
	\end{proof}
	

	\section{A strengthening of the Hurewicz property}

	T. Banakh and L. Zdomskyy \cite{Ba-Zd} considered the following property, stronger than the Hurewicz property: let $\E$ be a class of sets in a Polish space $Z$ containing all $G_\delta$-sets; $X\sub Z$ has {\sl the $\E$-separation property} if for every $E\in \E$ containing $X$ there is a $\s$-compact set containing $X$ and contained in $E$.

\ms 
    
	The following two results are related to this notion.
	
	\ms 
	
	The first one strengthens Theorem 2.1 in \cite{sz-w-z}. Let us recall that a subset $A$ of $\cantor$ is {\sl perfectly meager in the transitive sense} (PMT for short, cf. \cite{sz-w-z}) if,  for every perfect subset  $P$ of $\cantor$, there exists an $F_\sigma$-set $F$ in $\cantor$ such
	that $A \sub F$ and $F\cap (P+t)$ is meager in $P+t$ for each $t\in \cantor$. Using the natural identification of the product of $\cantor \times \cantor$ with $\cantor$, which preserves the algebraic and topological structure, we can apply this notion also to subsets of $\cantor \times \cantor$.
	
	\begin{theorem}\label{IMG, 30-12-2024}
		For every $\kappa$-Sierpiński set $S\sub \cantor$ of cardinality $\co$, there is $T\sub S$, $|T|=\co$, and a function $f: T\to \cantor$ such that, whenever $B$ is a Borel set in $\cantor\times\cantor$ containing the graph $Gr(f)$ of $f$, there is a $\s$-compact set $F$ with $|T\setminus F|<\kappa$ and $F\times \cantor\sub B$. In particular,  $Gr(f)$ is not PMT and  if $\kappa=\be$, $Gr(f)$ is a Hurewicz $\lambda'$-set in $\cantor\times\cantor$ which projects onto $\cantor$.		
	\end{theorem}

We start the proof with the following technical lemma.

\begin{lemma}\label{Hilg_for_full_measure}
Let $K$ be a compact subset of $\cantor$ of positive $\lambda$-measure, and let $T$ be a subset of $K$ intersecting each Borel (hence also each analytic) set in $K$ of positive $\lambda$-measure in $\co$ many points.

Then, for any Polish space $Y$, there exists a function $f: T\to Y$ such that, whenever $B$ is a Borel set in $\cantor\times Y$ containing 
$Gr(f)$,
there is a $\s$-compact set $F\subset K$ with $\lambda(F) = \lambda(K)$ and $F\times Y\sub B$.
\end{lemma}

\begin{proof}
We can split $T$  into  sets $T_\alpha$, $\alpha<\co$, intersecting each Borel set in $K$ of positive $\lambda$-measure.

Let $f:T\to Y$ be a Hilgers function associated with the set $A=K\times Y$, the family $\mathcal{F} = \{B_\alpha: \alpha < \co\}$  of all Borel sets in $\cantor\times Y$ such that for any $\alpha < \co$,  $\pi_{\cantor}(B_\alpha)$ contains $T$, and the partition  $\mathcal{P}= \{T_\alpha: \alpha < \co\}$, cf. Section \ref{Hilgers functions}. We shall check that $T$ and $f$ have the required properties.

Let $B$ be a Borel set in $\cantor\times Y$ containing  $Gr(f)$. 

 Then $T \sub \pi_{\cantor}(B)$, so there is $\alpha<\co$ such that $B=B_\alpha$ and, consequently, $Y\subseteq (B_\alpha)_x$  for all  $x\in T_\alpha$, i.e.,
 $T_\alpha \times Y \sub B_\alpha$ (cf. (\ref{Hilg_funct})).
  
Then $\pi_{\cantor}((K\times Y) \setminus B_\alpha)$ is an analytic subset of $K$ disjoint from $T_\alpha$, hence it is a $\lambda$-null set.  Let $N$ be a $\lambda$-null $G_\delta$-set containing  $\pi_{\cantor}((K\times Y) \setminus B_\alpha)$ and let 
$F = K \setminus N$. Then $F$ is \s-compact,  $\lambda(F) = \lambda(K)$, and  $F \times Y \sub B_\alpha=B$.
\end{proof}

\begin{proof}[Proof of Theorem \ref{IMG, 30-12-2024}]
Let $\A$ be a maximal pairwise disjoint collection of Borel sets of positive $\lambda$-measure in $\cantor$ intersecting $S$ in a set of cardinality less than $\co$. Then $\A$ is countable, so $|S\cap \bigcup\A|<\co$ and  $|S\setminus  \bigcup\A|=\co$.

Let $T^*=\cantor\setminus \bigcup\A$. Then $S\setminus \bigcup\A$ is a $\kappa$-Sierpiński set contained in $T^*$, hence $\lambda(T^*)>0$. Moreover, by the maximality of $\A$, each Borel set in $T^*$ of positive $\lambda$-measure intersects $S$ in $\co$ many points. In particular if we fix   a compact set $K\sub T^*$ with $\lambda(K)>0$ and let $T=K\cap S$, then $T\sub K$ is a $\kappa$-Sierpiński set of cardinality $\co$ intersecting each analytic set in $K$ of positive $\lambda$-measure in $\co$ many points.

Let $f:T\to \cantor$ be a  function given by Lemma \ref{Hilg_for_full_measure} applied for $K, T$ and $Y = \cantor$. We shall check that T and f have the required properties.

So let $B$ be a Borel set in $\cantor\times \cantor$ containing $Gr(f)$, and let $F\sub K$ be a $\s$-compact set as in Lemma \ref{Hilg_for_full_measure}. Since $\lambda(K\setminus F)=0$ and $T\subseteq K$ is a $\kappa$-Sierpiński set, $|T\setminus F|< \kappa$, as required. 

 To prove that  $Gr(f)$ is not PMT consider $P = \{\bar 0\} \times \cantor$ and let $B$ be any $F_\sigma$ (or even  Borel) set in $\cantor\times\cantor$ containing $Gr(f)$. Then, by what has already been proved,   there exists $a\in\cantor$ such that $\{a\} \times \cantor\sub B$ and   $B \cap ((a,{\bar 0})+P)=\{a\} \times \cantor$  is not meager in $(a,{\bar 0})+P=\{a\} \times \cantor$.

\ss

For the rest of the proof let us additionally assume that $\kappa=\be$. We keep the notation from the first part of the proof.

\ss 

 From $|T\setminus F|< \kappa=\be$ we have  $|Gr(f) \setminus (F \times \cantor)| < \be$, so  $|Gr(f) \setminus (F \times \cantor)|$ is a Hurewicz space. It follows that if $B$  is additionally assumed to be a $G_\delta$-set in $\cantor\times \cantor$ containing  $Gr(f)$, we can cover $Gr(f) \setminus (F \times \cantor)$ by a $\sigma$-compact set $L$ contained in $B$. Then $(F \times \cantor) \cup L$ is a $\sigma$-compact set containing $Gr(f)$ and contained in $B$ witnessing that $Gr(f)$ is a Hurewicz space.

\ss 
 
To prove that $Gr(f)$ is a $\lambda'$-set in $\cantor\times\cantor$, it is now enough by Lemma \ref{lambda'}, to show that it is a $\lambda$-space.  For this, let $A$ be an arbitrary countable subset of $Gr(f)$. Repeating the argument from the first part of the proof we pick a $\lambda$-null $G_\delta$-set  $N$ containing  $\pi(A)$ and letting $F = \cantor \setminus N$ we conclude that $F$ is \s-compact set disjoint from $\pi(A)$ with $|T\setminus F|< \be$. Consequently, letting $F_0 = F \times\cantor$ we have $F_0\sub (\cantor\times\cantor)\setminus A$ and
$|(Gr(f) \setminus (F_0 \cup A)| < \be$, so there is  an $F_\sigma$-set $F_1$ in $\cantor\times\cantor$ such that $(Gr(f) \setminus (F_0 \cup A) \sub F_1 \sub  (\cantor\times\cantor)\setminus A$. Finally, let $F= F_0 \cup F_1$. Then 
$Gr(f) \setminus A = F \cap Gr(f)$ which shows that $A$ is a  $G_\delta$-set  in $Gr(f)$.

\ss 
 Clearly, $Gr(f)$ projects onto  the second coordinate $\cantor$, since otherwise we would have $Gr(f) \sub \cantor \times (\cantor\setminus\{a\})$ for some $a\in \cantor$, and then $B =  \cantor \times  (\cantor\setminus\{a\})$ would be a Borel set containing  $Gr(f)$ and $B_x\neq\cantor$ for all $x$, contrary to the earlier conclusion.
\end{proof}

	\ms 
	
	The second result, where $\B$ ($\A$, respectively) stands for the family of Borel (analytic, respectively)  sets in $\cantor\times\cantor$, provides a solution (under $V=L$) of the second part of Problem 3 in \cite{Ba-Zd}.
	
	\begin{theorem}\label{IMG, 9-12-2024}
		Assuming $V=L$, there is a set in $\cantor$ with the $\B$-separation property but without the $\A$-separation property. 
		
	\end{theorem}

In the proof we shall use the following fact essentially due to G\"odel. 

Let $\pi:\cantor\times\cantor\to\cantor$ denote the projection onto the first axis.

\begin{lemma}\label{Godel}
Assuming $V=L$, there is a coanalytic set $C\sub \cantor\times\cantor$ such that
\begin{enumerate}[(i)]
\item $\pi(C)=\cantor$ and 
$|\pi^{-1}(t)\cap C|=1$ for each $t\in \cantor$,

\item for every analytic set $E\sub C$, $\pi(E)$ is countable.

\end{enumerate}

\end{lemma}

\begin{proof}
Let $W_0\sub \cantor$ be a ${\mathbf \Delta_2^1}$ Bernstein set (cf. \cite[Example 8.24]{ke}) which means that neither $W_0$ nor $W_1=\cantor\setminus W_0$  contains a non-empty perfect subset of $\cantor$ (the proof of the existence of such a set under $V=L$  is based on the work of G\"odel, cf. \cite{khom}).

By the Novikov-Kond\^o uniformization theorem (cf. \cite[Theorem 36.14]{ke}), there are disjoint coanalytic sets $C_i\sub\cantor\times\cantor$, $i=0,1$, with $\pi(C_i)=W_i$ and $|\pi^{-1}(t)\cap C_i|=1$ for each $t\in W_i$. Then $C=C_0\cup C_1$ has required properties. 

Indeed, if $E \sub C$ is an analytic set, then $E_i = E\cap C_i =
E \cap ((2^N \times 2^N) \setminus C_{i-1})$ is an analytic subset of $C_i$. Consequently,  $\pi(E\cap C_i)$ is an analytic subset of $W_i$. Since  $W_i$ is a Bernstein set, $\pi(E\cap C_i)$ is countable, and so is  $\pi(E\cap C)=\pi(E\cap C_1)\cup \pi(E\cap C_2)$.
\end{proof}

We  now   proceed to the proof of Theorem \ref{IMG, 9-12-2024}.

\begin{proof}
With the help of CH we first inductively construct a Sierpiński set $S$ in $\cantor$ which intersects every analytic non-null set in $\cantor$, then we partition it into  non-empty sets $S_\alpha$, $\alpha<\omega_1$, with the same property,
and finally we let $\mathcal{P} = \{S_\alpha: \alpha < \omega_1\}$.

Let $A=(\ctc) \setminus C$ where $C$ is the set given by Lemma \ref{Godel},  and let $\mathcal{F} = \{B_\alpha: \alpha < \omega_1\}$ be the family of all Borel sets in $\cantor\times \cantor$ such that, for any $\alpha < \omega_1$,  $\pi(B_\alpha)$ contains $S$.

Let $f: S \to \cantor$ be a Hilgers function associated with the set $A$, the family $\mathcal{F}$ and the partition  $\mathcal{P}$, cf. Section \ref{Hilgers functions}. We shall check that $Gr(f)$ has the required properties.

\ss 

To show that $Gr(f)$ has the $\B$-separation property let us consider a Borel set $B$ in $\cantor\times\cantor$ containing $Gr(f)$. Then there is $\alpha<\omega_1$ such that $B=B_\alpha$ and, consequently, $A_x\subseteq (B_\alpha)_x$  for all  $x\in S_\alpha$, cf. (H1) in subsection \ref{Hilgers functions}. Therefore, the projection $\pi(A\setminus B)$  is an analytic set ($A$ being analytic as the complement of $C$)  disjoint from $S_\alpha$, hence it is a $\lambda$-null set. 

Let $N_1$ be a $\lambda$-null $G_\delta$-set containing $\pi(A\setminus B)$ and let $F_1=\cantor\setminus N_1$. The set $F_1$ is \s-compact and $\pi^{-1}(F_1)\cap A\sub B$, which implies that the set $E=\pi^{-1}(F_1)\setminus B$ is contained in $(\cantor\times\cantor)\setminus A=C$. Since, clearly, $E$ is Borel, $\pi(E)$ is countable, by Lemma \ref{Godel}. 

Let $N_2$ be a $\lambda$-null  $G_\delta$-set in $\cantor$ containing $\pi(E)$ and let $F_2=F_1\setminus N_2$. Then the set $F_2$ is \s-compact, $\pi^{-1}(F_2)\sub B$ and, since $\lambda(F_2)=1$ and $S$ is a Sierpiński set, $|S\setminus F_2|\leq\aleph_0$. Consequently, $|Gr(f)\setminus \pi^{-1}(F_2)|\leq \aleph_0$, hence $F=\pi^{-1}(F_2)\cup (Gr(f)\setminus \pi^{-1}(F_2))$ is a $\sigma$-compact set containing $Gr(f)$ and contained in $B$ witnessing that $Gr(f)$ has the $\B$-separation property.

\ss 

The set $A$ containing $Gr(f)$ is analytic and $A_x\neq \cantor$ for all $x\in \cantor$. Since, as we have already established, for each Borel set $B$ containing $Gr(f)$ there is $x$ with $B_x=\cantor$ (in fact, $B_x=\cantor$ for every $x\in F_2$), there is no Borel set containing $Gr(f)$ and contained in $A$.

In particular, $Gr(f)$ fails the $\A$-separation property.
\end{proof}


	\section{Products of $C$-spaces}
	
Recall that a separable metrizable space $X$ is a {\sl $C$-space} if for every sequence $\G_1, \G_2,\ldots$ of open covers of $X$ there is a sequence $\H_1, \H_2,\ldots$ of families of pairwise disjoint open subsets of $X$ such that $\H_i$ refines $\G_i$ for $i=1,2,\ldots$  (i.e.,  each member of $\H_i$ is contained in a member of $\G_i$) and the union $\bigcup_{i=1}^\infty \H_i$ is a cover of $X$, cf. \cite{eng-2}. Note that, when we want to verify that a subspace $X$ of a metrizable space $Z$ is a $C$-space, we can require that the families $\G_i,\H_i$ in the above definition consists of sets open in $Z$. 
This follows from the well-known fact that, for each disjoint collection $\U$ of open subsets of a subspace $X$ of a metrizable space $Z$, there is a disjoint collection $\U'$ of open subsets of $Z$ such that $\U = \{U\cap X: U\in \U'\}$.
	
The main results of this section address some problems concerning the products of $C$-spaces  formulated by M. Sakai and M. Scheepers \cite{sak-sch}. The first one (Theorem  \ref{RP, 23-11-2024}) provides a negative answer to Problem 6.13 in \cite{sak-sch}, while the second one (Theorem  \ref{IMG, 27-11-2024}) is a negative solution of Problem 6.6 in \cite{sak-sch}.

	\begin{theorem}\label{RP, 23-11-2024}
		Assuming CH, for each Sierpiński set $S\sub \cantor$ there is a Hurewicz $C$-space $H\sub I^\NN$ such that $S\times H$ is not a $C$-space. Moreover, $H$ is a $\lambda$-space.
	\end{theorem}

We will use the following simple observation 
(recall that a subspace of a $C$-space need not be a $C$-space).

\begin{proposition}\label{Hur_sub_C}
Each Hurewicz subspace  $X$  of a compact metrizable $C$-space $K$ is a $C$-space.
\end{proposition}

	\begin{proof}
Let $\G_1, \G_2,\ldots$ be a sequence  of open in $K$ covers of $X$. Let $G = \bigcap\{\bigcup \G_i: i = 1,2,\dots\}$. Since $X$ is a Hurewicz space,  we can find a $\s$-compact set $F$ with $X\sub F \sub G$ (cf. subsection \ref{Hur_and_Meng}). The subspace $F$ is a $C$-space, cf.\ 
\cite[proof of 6.3.8]{eng-2}, therefore we can find a sequence $\H_1, \H_2,\ldots$ of families of pairwise disjoint open subsets of $F$ such that $\H_i$ refines $\G_i$ for $i=1,2,\ldots$ and the union $\bigcup_{i=1}^\infty \H_i$ covers $F$, hence also $X$.
	\end{proof}

We also need the next lemma.

\begin{lemma}\label{lem_C_space1}
Let $X,Y$ be separable metrizable spaces, and let $Z$ be a subspace of $X\times Y$ such that, for each open set $V$ in $X\times Y$ containing $Z$, there is $x\in X$ with $\{x\}\times Y \sub V$. 

If $Z$ is a $C$-space, then so is $Y$. 
\end{lemma}

\begin{proof}
Let $\W_1, \W_2,\ldots$ be a sequence of open covers of $Y$.

We define $\U_i = \{X\times W: W\in \W_i\}, i=1,2,\dots$. Using the $C$-space property of $Z$ we can find open in $X\times Y$ disjoint collections $\V_1, \V_2,\ldots$ such that $\V_i$ refines $\U_i$ and $Z \sub V = \bigcup\{\bigcup\V_i: i=1,2,\dots\}$. By our assumption on $Z$ there is $x\in X$ with $\{x\}\times Y \sub V$. Then $\D_i = \{V_x: V\in \V_i\}$ is an open disjoint refinement of $\W_i$ and $Y\sub \bigcup\{\bigcup\D_i: i=1,2,\dots\}$.
\end{proof}

	\begin{proof}[Proof of Theorem \ref{RP, 23-11-2024}]
 There exists a metrizable compact  $C$-space $L$ containing a $G_\delta$-subspace $G$ which is not a $C$-space (
 cf.\ \cite[Example 6.1.21, Theorem 6.3.10 and Problem 6.3.D(b)]{eng-2}).

Let $\pi: \cantor\times L \to \cantor$ be the projection.

Let $T\sub \cantor$ be a Sierpiński set disjoint from $S$, with $\lambda^*(T) = 1$, given by Lemma \ref{Sierp full measure}. Observe that, by CH, $T$ intersects each Borel set in $\cantor$ of positive $\lambda$-measure in $\co$ many points. Therefore, we can apply Lemma \ref{Hilg_for_full_measure}  for $T$, $K = \cantor$ and $Y = L$, obtaining a function $f:T\to L$ such that
\begin{eqnarray}\label{5.1.1}
&&\mbox{for each Borel set } B\sub \cantor\times L \mbox{ containing }  Gr(f) \mbox{ there is} \\
\nonumber &&\mbox{a $\s$-compact set } F\sub \cantor \mbox{ with } \lambda(F) =1 \mbox{ and } F\times L \sub B.
\end{eqnarray}
Let $g:S\to G$ be a Hilgers function associated with the set $S\times G$ and  the family $\mathcal{F} = \{V_\alpha: \alpha < \co\}$  of all open sets in $\cantor\times L$ such that for any $\alpha < \co$,  $\pi(V_\alpha)$ contains $S$. By condition (\ref{Hilg_funct}), cf.\ subsection \ref{Hilgers functions}, 
$Gr(g)$
has the the following property
\begin{eqnarray}\label{5.1.2}
&&\mbox{for each open } V\sub \cantor\times L \mbox{ containing }  Gr(g) \mbox{ there is} \\
\nonumber &&s\in S   \mbox{ such that } \{s\}\times G \sub V.
\end{eqnarray}

Finally, let
\begin{equation}\label{5.1.3}
H = Gr(f) \cup Gr(g) \sub \cantor\times L.
\end{equation}

We claim that $H$ has the desired properties.

\ss  

First, we shall verify that $H$ is a Hurewicz space. So, take a $G_\delta$-subset $B$ of $\cantor\times L$ containing $H$. By (\ref{5.1.1}) we have a  $\s$-compact set $F\sub \cantor$  such that  $\lambda(F) =1 $ and  $F\times L \sub B$. The  set $E = (S\cup T)\setminus F$ is countable,  $S\cup T$ being a Sierpiński set. Therefore, the set 
$(F\times L)\cup[(E \times L)\cap H]$ is $\s$-compact, and we have
\begin{equation}\label{5.1.4}
H \sub (F\times L)\cup[(E \times L)\cap H] \sub B.
\end{equation}

Since the product $\cantor\times L$ is a $C$-space 
(cf.\ \cite[Theorem\ 6.3.11]{eng-2}),
Proposition \ref{Hur_sub_C} assures us that so is $H$.

Next, note that condition (\ref{5.1.2}) together with Lemma \ref{lem_C_space1} implies that
\begin{equation}\label{5.1.5}
\mbox{the graph } Gr(g) \mbox{ is not a  $C$-space.} 
\end{equation}

Finally, observe that the product $S\times H$ contains a closed copy $\{(s,(x,y))\in S\times H: s= x\}$ of $Gr(g)$, hence, by (\ref{5.1.5})  $S\times H$ is not a $C$-space.

Moreover, since the projection onto the first axis maps $H$ injectively into the $\lambda$-set $S\cup T$ ($S\cup T$ is a Sierpiński set, cf.\ Lemma \ref{Sierp_lambda'}), $H$ is a $\lambda$-set.
	\end{proof}

	\begin{theorem}\label{IMG, 27-11-2024}
		Assuming $\be=\co$, there is $E\sub I^\NN$ such that
		\begin{enumerate}[(i)]
			\item $E$ is a $C$-space but $E\times \NN^\NN$ fails this property,
			
			\item for each Hurewicz space $Y\sub I^\NN$, $E\times Y$ is Hurewicz; in particular, $E^n$ is Hurewicz for $n=1,2,\ldots$.
		\end{enumerate}
		
	\end{theorem}

For the proof we will need two auxiliary facts. The first one is a simple observation extracted from the proof of properties of Example 2 in \cite{ep}. Recall that, for a cardinal number $\kappa$, a space $X$ is {\sl $\kappa$-concentrated about a set $A\sub X$}, if, for each open set in $X$ containing $A$, we have $|X\setminus U| < \kappa$.

\begin{proposition}\label{lem_C_space2}
Let $X$ be a separable metrizable space $\co$-concentrated about a subspace $Y\sub X$. If $Y$ is a $C$-space, then so is $X$. 
\end{proposition}

\begin{proof}
Let $\W_1, \W_2,\ldots$ be a sequence of open covers of $X$. 

By the $C$-space property of $Y$, we can find open in $X$ disjoint collections $\V_2, \V_3,\ldots$ such that $\V_i$ refines $\W_i,\ i=2,3,\dots$ and $Y \sub V = \bigcup\{\bigcup\V_i: i=2,3,\dots\}$. Since $|X\setminus V| < \co$, the subspace $X\setminus V$ is zero-dimensional, hence there is an open in $X$ disjoint collection $\V_1$ refining $\W_1$ and covering $X\setminus V$. Clearly, the sequence $\V_1,\V_2, \V_3,\ldots$ witnesses the $C$-space property of $X$.
\end{proof}

In the next statement  (probably a part of the folklore in this area), $\mathfrak{B}$ is a $\be$-scale set described in subsection \ref{bscale}. This statement implies in particular, that for any Hurewicz space $Y\sub I^\NN$, the product $\mathfrak{B}\times Y$ is Hurewicz, which was established by Miller, Tsaban and Zdomskyy in \cite[Theorem 6.7]{Mi-Ts-Zd}.

\begin{lemma}\label{lem_Hur_prod}
Let $K$ be a metrizable compact space, $f:K\to \cantor$ be a continuous map, and $E=A\cup C$ be a subspace of $f^{-1}(\mathfrak{B})$ such that $C = f^{-1}(\bar{Q}), A\cap C = \emptyset$, and $|A\cap f^{-1}(t)| = 1$ for all $t\in \mathfrak{B}\setminus \bar{Q}$.

Then, for every Hurewicz space $Y\sub I^\NN$, the product $E\times Y$ is Hurewicz. 
\end{lemma}

\begin{proof}
Let $G\sub K\times I^\NN$ be a $G_\delta$-set containing $E\times Y$.

Let $\pi_1: K\times I^\NN \to K, \pi_2: K\times I^\NN \to I^\NN$ be the projections.

The set $L' = \pi_2((C\times I^\NN)\setminus G)$ is $\s$-compact and disjoint from $Y$. Since $Y$ is Hurewicz, there is a $\s$-compact set $L\sub I^\NN$ such that 
\begin{equation}\label{5.6.1}
Y\sub L \sub I^\NN \setminus L'.
\end{equation}
Then 
\begin{equation}\label{5.6.2}
C\times Y \sub C\times L \sub G.
\end{equation}
Therefore, the set
\begin{equation}\label{5.6.3}
M = \pi_1((K\times L)\setminus G) \mbox{ is $\s$-compact, and } M\cap C =\emptyset\,.
\end{equation}
Since $f(M)$ is a $\s$-compact subset of $\bar{P}$, by 
property (BS3) of $\mathfrak{B}$ (cf.\ subsection \ref{bscale}), there is a  $\s$-compact set $N$ in $\cantor$ with
\begin{equation}\label{5.6.4}
\bar{Q} \sub N \sub \cantor \setminus f(M) \mbox{\ and\ } |\mathfrak{B}\setminus N| < \be.
\end{equation}
Letting $F_1 = f^{-1}(N)\times L$, by (\ref{5.6.1}), (\ref{5.6.3}) and  (\ref{5.6.4}) we have
\begin{equation}\label{5.6.5}
C\times Y \sub f^{-1}(N)\times L\sub F_1\sub G \mbox{ and } |E\setminus f^{-1}(N)| < \be.
\end{equation}
Therefore, the product $(E\setminus f^{-1}(N))\times Y$ is Hurewicz, as the union of less than $\be$ Hurewicz sets of the form $ \{x\}\times Y$, where $x\in E\setminus f^{-1}(N)$ (cf. subsection \ref{Hur_and_Meng}). Hence, we can find a $\s$-compact set $F_2$ with
\begin{equation}\label{5.6.6}
(E\setminus f^{-1}(N))\times Y \sub F_2 \sub G.
\end{equation}
Finally, $F_1\cup F_2$ is a $\s$-compact set containing $E\times Y$ and contained in $G$.
\end{proof}

	\begin{proof}[Proof of Theorem \ref{IMG, 27-11-2024}]
First, we will recall a (slightly modified) construction from
\cite[Example 2]{ep}. 

Let $L = \cantor \times I^\NN$, and let $\pi_0: L\to \cantor,\ \pi_n:L \to I^n, n= 1,2,\dots$, be the projections defined by
\begin{eqnarray}\label{5.4.1}
\pi_0((x,(t_i)_{i=1}^\infty)) &=& x\,;\\
\pi_n((x,(t_i)_{i=1}^\infty)) &=& (t_1,t_2,\dots,t_n)  \,.\label{5.4.2}
\end{eqnarray}
for $x\in \cantor, (t_i)_{i=1}^\infty \in I^\NN, n= 1,2,\dots$. 

Fix a bijective enumeration $q_1,q_2,\dots$ of the set $\bar{Q}$.  
We consider the closed equivalence relation on the space $L$ which identifies points $(x,(t_i)_{i=1}^\infty)$ and  $(x',(t'_i)_{i=1}^\infty)$,  whenever there is $n\in \{1,2,\dots\}$ such that $x=x'=q_n$ and $\pi_n((q_n,(t_i)_{i=1}^\infty)) = \pi_n((q_n,(t'_i)_{i=1}^\infty)$. Denote the quotient compact space obtained in this way by $K$, and let $q:L\to K$ be the natural quotient map. Let $f:K\to \cantor$ be the (unique) continuous map such that $f\circ q = \pi_0$. 


Observe that 
\begin{eqnarray}\label{5.4.3}
&&q|(\bar{P}\times I^\NN) \mbox{ is a homeomorphic embedding onto } f^{-1}(\bar{P})\,;\\
&&f^{-1}(q_n) \mbox{ is  homeomorphic with } I^n, \mbox{ for } n= 1,2,\dots\,\label{5.4.4}
\end{eqnarray}
Let $g:\mathfrak{B}\setminus \bar{Q} \to I^\NN$ be a Hilgers function associated with the set $(\mathfrak{B}\setminus \bar{Q}) \times I^\NN$ and  the family $\mathcal{F} = \{V_\alpha: \alpha < \co\}$  of all open sets in $\bar{P}\times I^\NN$ such that for any $\alpha < \co$,  $\pi_0(V_\alpha)$ contains $\mathfrak{B}\setminus \bar{Q}$. By condition (\ref{Hilg_funct}), cf.\ subsection \ref{Hilgers functions}, 
$Gr(g)$
has the the following property
\begin{eqnarray}\label{5.4.5}
&&\mbox{for each open } V\sub \bar{P}\times I^\NN \mbox{ containing }  Gr(g) \mbox{ there is} \\
\nonumber &&x\in  \mathfrak{B}\setminus \bar{Q}  \mbox{ such that } \{x\}\times I^\NN \sub V.
\end{eqnarray}
Now, we can define our space $E$ in the following way
\begin{equation}\label{5.4.6}
E = A\cup C \sub K,\quad \mbox{where } A = q(Gr(g))\mbox{ and } C = f^{-1}(\bar{Q})\,.
\end{equation}
Clearly, 
\begin{equation}\label{5.4.7}
|A\cap f^{-1}(x)| = 1 \quad \mbox{for all } x\in \mathfrak{B}\setminus \bar{Q}\,.
\end{equation}

Observe that the space $E$ is $\be$-concentrated about the set $C$. Indeed, if $U$ is an open set in $K$ containing $C$, then $f(K\setminus U)$ is a compact subset of $\cantor$ disjoint from $\bar{Q}$, hence, by property (BS2) of $\mathfrak{B}$ (cf.\ subsection \ref{bscale}), $|f(K\setminus U)\cap \mathfrak{B}|< \be$. Consequently, $|E\setminus U| < \be$, since $f$ is one-to-one on $A$, cf.\ (\ref{5.4.7}).

The set $C$ is a $C$-space, being a countable union of finite-dimensional compacta, cf.\ (\ref{5.4.4}), (\ref{5.4.6}) and \cite[Theorem 6.3.8]{eng-2}. From Proposition \ref{lem_C_space2} we conclude that $E$ is a $C$-space.

Since the Hilbert cube $I^\NN$ is not a $C$-space, Lemma \ref{lem_C_space1} and condition (\ref{5.4.5}) imply that the graph $Gr(g)$ is not a $C$-space, and neither is $A$, cf.\ (\ref{5.4.3}). The space $A$ homeomorphically embeds in $E\times \bar{P}$ as the closed subspace $\{(x,f(x)):\ x\in A\}=\{(x,y)\in E\times \bar{P}:\ y=f(x)\}$, and since $\bar{P}$ is homeomorphic to $\baire$, we get (i) of the theorem.

Statement (ii) of the theorem follows immediately from Lemma \ref{lem_Hur_prod}, cf.\ conditions (\ref{5.4.6}), (\ref{5.4.7}).
	\end{proof}




\begin{thebibliography}{JMSS}

		
		\bibitem[BS]{b-s} T.~Bartoszyński, S.~Shelah, {\it Continuous images of sets of reals},   Topology  Appl.  \textbf{116(2)}  (2001),  243--253. 	
		
		
		
		
		
		
		
		
		
		
		
		
		
		
		\bibitem[BZ]{Ba-Zd} T.~Banakh, L.~Zdomskyy, \textit{Separation properties between the $\s$-compactness and Hurewicz property}, Topology Appl. \textbf{156} (2008), 10--15.

        \bibitem[Bla]{Bla}
A.~Blass, \textit{Combinatorial cardinal characteristics of the continuum} in: {\it Handbook of Set Theory} (M. Foreman, A. Kanamori, and M. Magidor, eds.),
Springer, 2010, 395--491.
		
		
		
		
		
		
		\bibitem[Eng2]{eng-2} R.~Engelking, \textit{Theory od dimensions. Finite and infinite},   Heldermann Verlag, Berlin, 1989. 
		
		\bibitem[Eng1]{eng-1} R.~Engelking, \textit{General Topology}, Sigma Series in Pure Mathematics,  Heldermann Verlag, 1989. 
		
		\bibitem[Hi]{hil} A.~Hilgers, \textit{Bemerkung Zur Dimensionstheorie}, Fund. Math. \textbf{28} (1937), 303–-304.

       
		
		
		
		
		
		
		
		
		
		\bibitem[Hu1]{hur-1} W.~Hurewicz, \textit{\" Uber  eine Verallgemeinerung des Borelschen Theorems}, Math. Z. 24 (1925), 401–-421.
		
		 \bibitem[Hu2]{hur-2} W.~Hurewicz, \textit{\" Uber Folgen stetiger Funktionen}, Fund. Math. \textbf{9} (1927), 193–-204.
		
		\bibitem[JMSS]{j-m-s-s} W. Just, A.W. Miller, M. Scheepers, P.J. Szeptycki, \textit{The combinatorics of open covers (II)}, Topology Appl. \textbf{73} (1996), 241--266.
		
		
		
		\bibitem[Ke]{ke} A.~S.~Kechris, \textit{Classical descriptive set theory}, Graduate Texts in Math. 156, Springer-Verlag, 1995.

        
        \bibitem[Kh]{khom} Y.~Khomskii, \textit{Regularity Properties and Definability
in the Real Number Continuum}, https://www.math.uni-hamburg.de/home/loewe/pdf/khomskii.yurii.pdf.

       \bibitem[Me]{meng} K.~Menger, \textit{Einige \" Uberdeckungss\" atze der Punktmengenlehre}, Sitzungsber. Wien. Akad. 133 (1924) 421–-444.
		
		
		
      \bibitem[Mich]{mich} E.~Michael, \textit{Paracompactness and the Lindel\" of property in finite and countable Cartesian products}, Comp. Math. \textbf{23(2)} (1971), 199--214.

        
		\bibitem[Mi]{mi-1} A.~W.~Miller, \textit{Special subsets of the real line} in \textit{Handbook of set--theoretic topology}, North-Holland 1984, 201--233.
		
		
		
		
		
		
		
		
		
		
		
		
		\bibitem[MTZ]{Mi-Ts-Zd} A.~W.~Miller, B.~Tsaban, L.~Zdomskyy, \textit{Selective covering properties of product spaces},	Ann. Pure Appl. Log.  \textbf{165} (2014), 1034–-1057.
		
		
		
	
		
		\bibitem[EPo]{ep} E.~Pol, \textit{A weakly infinite-dimensional space whose product with the irrationals is strongly infinite-dimensional}, Proc. Amer. Math. Soc. \textbf{98} (1986), 349--352.

		
		
		\bibitem[PP]{p-p} E.~Pol, R.~Pol, \textit{On metric spaces with the Haver property which are Menger spaces},  Topology Appl.    \textbf{157} (2010), 1495-–1505.
		
		\bibitem[PZ]{p-z} R.~Pol, P.~Zakrzewski, \textit{Countably perfectly meager sets}, J. Symbolic Logic  \textbf{86(3)} (2021), 1--17.

        	
        

		
		\bibitem[SS]{sak-sch} M.~Sakai, M.~Scheepers, \textit{The combinatorics of open covers}, in \textit{Recent progress in General Topology III}, Atlantis Press 2014.
		
		
		
		
		\bibitem[ST1]{sz-ts-1} P.~Szewczak, B.~Tsaban, {\it Products of Menger spaces: a combinatorial approach},  Ann. Pure Appl. Log.  \textbf{168} (2017), 1--18.
		
		\bibitem[ST2]{sz-ts-2} P.~Szewczak, B.~Tsaban, {\it Products of general Menger spaces},   Topology Appl. \textbf{255} (2019), 41--55.  .
		
		\bibitem[SWZ]{sz-w-z} P.~Szewczak, T.~Weiss, L.~Zdomskyy, {\it Small Hurewicz and Menger sets which have large continuous images},  https://arxiv.org/abs/2406.12609.


        \bibitem[Ts]{ts} B.~Tsaban, {\it Menger’s and Hurewicz’s Problems:  Solutions from “The Book” and refinements},   Contemp. Math. \textbf{533} (2011), 211--226.
		
		
		
		
		
		
		
	\end{thebibliography}
\end{document}